\documentclass[preprint,12pt]{elsarticle}\usepackage{amsmath,amssymb}
\newtheorem{defin}{Definition}[section]\newtheorem{thm}[defin]{Proposition}\newtheorem{corollary}[defin]{Corollary}\newproof{pf}{Proof}\newtheorem{lemma}[defin]{Lemma}

\begin{document}
\begin{frontmatter}

\title{\bf A differential calculus on Z$_3$-graded quantum superspace ${\mathbb R}_q(2|1)$}
\date{\today}

\author[sacelik]{Salih Celik\corref{cor}}\ead{sacelik@yildiz.edu.tr}
\cortext[cor]{Corresponding Author}

\address[sacelik]{Department of Mathematics,Yildiz Technical University, DAVUTPASA-Esenler
Istanbul, 34210 TURKEY.}

\begin{abstract}
We introduce a Z$_3$-graded quantum $(2+1)$-superspace and define Z$_3$-graded Hopf algebra structure on algebra of functions on the Z$_3$-graded quantum superspace. We construct a differential calculus on the Z$_3$-graded quantum superspace, and obtain the corresponding Z$_3$-graded Lie superalgebra.
We also find a new Z$_3$-graded quantum supergroup which is a symmetry group of this calculus.
\end{abstract}

\end{frontmatter}

\section{Introduction}\label{sec1}

Noncommutative geometry \cite{Con}, in recent years, started to play an important role in many different areas of mathematics and mathematical physics. The basic structure leading to the noncommutative geometry is a differential calculus on an associative algebra. The noncommutative differential geometry of quantum groups was introduced by Woronowicz \cite{Woro}. In his approach the differential calculus on the group is deduced from the properties of the group and it involves functions on the group, differentials, differential forms and derivatives.

After presentation of the quantum plane by Manin \cite{Manin1}, Wess and Zumino \cite{WZ} has developed a differential calculus on the quantum (hyper)plane which is covariant under the action of quantum group. In that approach, differential forms are defined in terms of noncommuting coordinates, and the differential and algebraic properties of quantum groups acting on these spaces are obtained from the properties of the spaces. The natural extension of their scheme to superspace \cite{Manin2} was introduced in e.g. \cite{Soni} and \cite{Celik1}.

Recently, there have been many attempts to generalize Z$_2$-graded constructions to the Z$_3$-graded case \cite{Kerner1}, \cite{Chung}, \cite{Kerner2}, \cite{Celik2}, \cite{AbBa}, \cite{Celik3} and Z$_N$-graded case \cite{Dub} and \cite{Nie}. The noncommutative geometry of the Z$_3$-graded quantum superplane is investigated in \cite{Celik2}. In this work, we will set up a differential calculus on the Z$_3$-graded quantum $(2+1)$-superspace, denoted by ${\mathbb R}_q(2|1)$. The differential calculus on the Z$_3$-graded quantum superspace involves functions on the superspace, differentials, differential forms and derivatives. The most important property of this calculus is that the operator {\sf d} satisfies ${\sf d}^3=0$ \, (${\sf d}^2\ne0$) and it contains as a consequence, not only first differentials ${\sf d}x^i$, $i = 1,2,3$, but involves also higher order differentials ${\sf d}^k x^i$, $k=1,2$. This paper is organized as follows.
In section 2, we define Z$_3$-graded quantum superspace and we introduce its Hopf algebra structure. In section 3, we give a differential calculus on the Z$_3$-graded quantum superspace. Z$_3$-graded quantum Lie superalgebra is given in section 4. The dual Hopf algebra is given in section 5.
In section 6, we introduce a new Z$_3$-graded quantum supergroup which is the symmetry group of the existing calculus. In section 7, we give an $R$-matrix which is associated with the quantum Z$_3$-graded superspace ${\mathbb R}_q(2|1)$.

\section{Z$_3$-graded quantum superspaces}\label{sec2}

In this section, we will state some properties of the Z$_3$-graded quantum superspaces.

\subsection{Z$_3$ gradation}\label{sec2.1}

The cyclic group Z$_3$ can be represented in the complex plane as multiplication by primary cubic roots of unity $j=e^{\frac{2\pi i}{3}}$ \, $(i=\sqrt{-1})$
$$1+j+j^2=0 \quad \mbox{and} \quad j^3=1.$$
In the following sections, we will continue taking $j=q$.

Let us consider a vector space $V$ of finite dimension $d$ with the direct sum decomposition $V=V_0\oplus V_1\oplus V_2$ such that there is a function with values in the group Z$_3$ given by $p(x)\equiv\alpha$ (mod 3) when $x\in V_\alpha$. This direct sum defines a Z$_3$-gradation on $V$ and any vector $v\in V$ is uniquely determined by the sum $v=v_0\oplus v_1\oplus v_2$ with its homogeneous components $v_\alpha \in V$.

\subsection{The Algebra of functions on the Z$_3$-graded quantum superspace}\label{sec2.2}

One can define a Z$_2$-graded superspace ${\mathbb R}(2|1)$ by dividing the superspace ${\mathbb R}(2|1)$ of 3x1 vectors into two parts ${\mathbb R}(2|1)=V_0\oplus V_1$. A vector is an element of $V_0$ (resp. $V_1$) and is of grade 0 (resp. 1) if it has the form given below:
\begin{equation}
\begin{pmatrix} x \\ y \\ 0\end{pmatrix}, \quad \mbox{resp.} \quad \begin{pmatrix} 0 \\ 0 \\ \theta \end{pmatrix}. \nonumber\end{equation}

We know that the Z$_2$-graded quantum superspace ${\mathbb R}_q(2|1)$ is defined as an associative, unital algebra generated by two even coordinates $x$, $y$ and an odd (Grassmann) coordinate $\theta$ with four quadratic relations \cite{Manin2}
\begin{equation*}
xy-qyx=0, \quad x\theta-q\theta x=0, \quad y\theta-q\theta y=0, \quad \theta^2=0
\end{equation*}
where $q$ is a nonzero complex number. As an algebra, $O({\mathbb R}_q(2|1))$ is the quotient of the super $q$-commutative polynomial algebra ${\mathbb C}[x,y,\theta]$ in three indeterminates $x$, $y$, $\theta$ (they are the coordinate functions on ${\mathbb R}_q(2|1)$) by the two-sided ideal generated by the elements $xy-qyx$, $x\theta-q\theta x$, $y\theta-q\theta y$, $\theta^2$. In the limit $q\to1$, this algebra is super-commutative and can be considered as the algebra of polynomials ${\mathbb R}(2|1)$ over the usual superspace, where $x$, $y$ and $\theta$ are the three coordinate functions.

One of possible ways to generalize the Z$_2$-graded superspace is given by the following definition.

\begin{defin}
Let $K\{x,y,\theta\}$ is a free associative algebra generated by $x$, $y$, $\theta$ and $I$ is a two-sided ideal generated by $xy-qyx$, $x\theta-q\theta x$, $y\theta-q^2\theta y$, $y^3$ and $\theta^3$. The quantum superspace ${\mathbb R}_q(2|1)$ with the function algebra
$$O({\mathbb R}_q(2|1)) = K\{x,y,\theta\}/I$$
is called Z$_3$-graded quantum superspace.
\end{defin}
In the Z$_3$-graded $O({\mathbb R}_q(2|1))$ the generator $x$ is of degree zero, the generator $y$ is of degree 1 and the generator $\theta$ is of degree 2.

In accordance with Definition 2.1, we have
\begin{equation} \label{1}
{\mathbb R}_q(2|1) = \left\{\begin{pmatrix} x \\ y \\ \theta \end{pmatrix}: xy=qyx, \,\,\, x\theta=q\theta x, \,\,\, y\theta=q^2\theta y, \,\,\, y^3=0=\theta^3\right\}\end{equation}
where $q$ is a cubic root of unity $(q^3=1)$.
The associative algebra $O({\mathbb R}_q(2|1))$ is $q$-commutative and known as the algebra of polynomials over the Z$_3$-graded superspace, where $x$, $y$ and $\theta$ are the three coordinate functions.

So, we define a Z$_3$-graded quantum superspace ${\mathbb R}_q(2|1)$ by dividing the quantum superspace ${\mathbb R}_q(2|1)$ of 3x1 vectors into three parts ${\mathbb R}(2|1)=V_0\oplus V_1\oplus V_2$. A vector is an element of $V_0$ (resp. $V_1$, $V_2$) and is of grade 0 (resp. 1, 2) if it has the form given below:
\begin{equation}
\begin{pmatrix} x \\0 \\ 0\end{pmatrix}, \quad \mbox{resp.} \quad \begin{pmatrix} 0 \\ y \\ 0 \end{pmatrix}, \quad \begin{pmatrix} 0 \\ 0 \\ \theta\end{pmatrix}. \nonumber\end{equation}

\noindent{\bf Note 1.} It is easy seen, from (\ref{1}), that the quantity $x^3$ is a central element of the algebra $O({\mathbb R}_q(2|1))$. As we will see in the next section, the cube of $x$ is in fact a central element for all differential algebra.

We also need to consider a quantized version of a dual algebra in three variables:
\begin{defin}
The quantum superspace ${\mathbb R}^*_q(2|1)$ with the function algebra
\begin{equation}
O({\mathbb R}^*_q(2|1)) = K\{\xi,\eta,z\}/(\xi\eta-q^2\eta\xi, \quad \xi z-z\xi, \quad \eta z-z\eta, \quad \xi^3, \quad \eta^3)\nonumber
\end{equation}
is called dual Z$_3$-graded quantum superspace.
\end{defin}
Hence, in accordance with Definition 2.2, we have
\begin{eqnarray} \label{2}
{\mathbb R}^*_q(2|1)=\left\{\begin{pmatrix} \xi \\ \eta \\ z\end{pmatrix}: \xi\eta=q^2\eta\xi, \quad \xi z=z\xi, \quad \eta z=z\eta, \quad \xi^3=0=\eta^3\right\}.
\end{eqnarray}
Here, the generators $\xi$, $\eta$, $z$ are of degree 1, 2, 0, respectively.

\noindent{\bf Note 2.} It is easy seen, from (\ref{2}), that the generator $z$ is a central element of the dual algebra $O({\mathbb R}^*_q(2|1))$.

In the next section, while we set up a differential calculus, the relations in (\ref{2}) will appear in Proposition 3.4.

\subsection{The Hopf algebra structure on ${\cal A}$}\label{sec2.3}

We define the extended Z$_3$-graded quantum superspace to be the algebra containing ${\mathbb R}_q(2|1)$, the unit and the inverse of $x$, $x^{-1}$, which obeys $xx^{-1}={\bf 1}=x^{-1}x$. We will denote the unital extension of $O({\mathbb R}_q(2|1))$ by ${\cal A}$. Let us begin with the following definition given by Majid \cite{Majid}.

\begin{defin}
For the Z$_3$-graded algebra ${\cal A}$, the product rule in the Z$_3$-graded algebra ${\cal A}\otimes{\cal A}$ is defined by
\begin{equation} \label{3}
(a\otimes b)(c\otimes d)=q^{p(b)p(c)}(ac\otimes bd)
\end{equation}
where $a,b,c,d$ are homogeneous elements in the algebra ${\cal A}$ whose grading is given by the function $p(.)$ with values in the set Z$_3=\{0,1,2\}$.
\end{defin}

\begin{thm}
The algebra ${\cal A}$ is a Z$_3$-graded Hopf algebra. The definitions of a coproduct, a counit and a coinverse on the algebra A are as follows:
{\bf (1)} The coproduct
$\Delta: {\cal A} \longrightarrow {\cal A} \otimes {\cal A}$
is defined by
\begin{equation} \label{4}
\Delta(x) = x \otimes x, \quad \Delta(y) = y \otimes {\bf 1} + x\otimes y, \quad \Delta(\theta) = \theta \otimes x + x^{-1} \otimes \theta.
\end{equation}

\noindent{\bf (2)} The counit $\epsilon: {\cal A} \longrightarrow {\mathbb C}$
is given by
\begin{equation} \label{5}
\epsilon(x)=1, \quad \epsilon(y)=0, \quad \epsilon(\theta)=0.
\end{equation}

\noindent{\bf (3)} If we extend the algebra $\cal A$ by adding the inverse of $x$ then the algebra $\cal A$ admits a ${\mathbb C}$-algebra antihomomorphism
(coinverse)
$\kappa: {\cal A} \longrightarrow {\cal A}$ defined by
\begin{equation} \label{6}
\kappa(x) = x^{-1}, \quad \kappa(y) = - x^{-1} y, \quad \kappa(\theta) = -q\theta.
\end{equation}
\end{thm}

\begin{pf}
It is not difficult to verify the following properties of the costructures:
The coproduct
$\Delta$ is coassociative in the sense that
\begin{equation} \label{7}
(\Delta \otimes \textrm{id}) \circ \Delta = (\textrm{id} \otimes \Delta) \circ \Delta
\end{equation}
where id denotes the identity map on ${\cal A}$ and $\Delta(ab)=\Delta(a) \Delta(b)$, $\Delta({\bf 1})={\bf 1}\otimes {\bf 1}$.

The counit $\epsilon$ has the property
\begin{equation} \label{8}
m \circ (\epsilon \otimes \mbox{id}) \circ \Delta = \mbox{id} = m \circ (\mbox{id} \otimes \epsilon) \circ \Delta
\end{equation}
where $m$ stands for the algebra product ${\cal A} \otimes {\cal A} \longrightarrow {\cal A}$ and $\epsilon(ab)=\epsilon(a) \epsilon(b)$, $\epsilon({\bf 1})=1$.
The coinverse $\kappa$ satisfies
\begin{equation} \label{9}
m \circ (\kappa \otimes \mbox{id}) \circ \Delta = \epsilon = m \circ (\mbox{id} \otimes \kappa) \circ \Delta
\end{equation}
and $\kappa(ab)=\kappa(b) \kappa(a)$, $\kappa({\bf 1})={\bf 1}$.
\end{pf}

The set $\{x^ky^l\theta^m: \, k\in{\mathbb N}_0, \, l,m=0,1,2\}$ forms a vector space basis of $O({\mathbb R}_q(2|1))$. The formula (\ref{4}) gives the action of the coproduct $\Delta$ only on the generators. The action of $\Delta$ on product of generators can be calculated by taking into account that $\Delta$ is an algebra homomorphism.

\section{A differential calculus on the Z$_3$-graded quantum superspace}\label{sec3}

In this section, we will set up a differential calculus on the Z$_3$-graded quantum superspace ${\mathbb R}_q(2|1)$. This calculus involves functions on this superspace, differentials and differential forms. It is sufficient to define the action of the differential operator on the coordinates and on their products.

\subsection{Differential algebra on superspace} \label{sec3.1}

We begin with the definition of the Z$_3$-graded differential calculus. Let $p(\alpha)$ denote the grade of $\alpha$.
\begin{defin}
Let ${\cal A}$ be an arbitrary associative (in general, noncommutative) algebra and $\Gamma^{\wedge n}$ be a space of $n$-form and ${\cal A}$-bimodule.
A Z$_3$-graded differential calculus on the algebra ${\cal A}$ is a Z$_3$-graded algebra $\Gamma^\wedge=\bigoplus_{n=0}^\infty \Gamma^{\wedge n}$ with a ${\mathbb C}$ linear exterior differential operator {\sf d} which defines the map ${\sf d}:\Gamma^\wedge \longrightarrow \Gamma^\wedge$ of grade one. A generalization of a usual differential calculus leads to the rules:
\begin{eqnarray} \label{10}
{\sf d}^3 & =& 0, \qquad ({\sf d}^2\ne0) \nonumber\\
{\sf d}(\alpha\wedge\beta) &=& ({\sf d}\alpha)\wedge\beta + q^{p(\alpha)} \, \alpha\wedge({\sf d}\beta), \\
{\sf d}^2(\alpha\wedge\beta) &=& ({\sf d}^2\alpha)\wedge\beta + (q^{p(\alpha)}+q^{p(d\alpha)}) \, (d\alpha)\wedge({\sf d}\beta) + q^{2p(\alpha)} \, \alpha\wedge({\sf d}^2\beta)\nonumber
\end{eqnarray}
for $\alpha\in\Gamma^{\wedge n}$ and $\beta\in\Gamma^{\wedge}$.
\end{defin}
For the product by elements $a\in \Gamma^{\wedge0}={\cal A}$ we shall write simply $a\rho$ and $\rho a$, $\rho\in \Gamma^{\wedge}$. The second condition in Definition 3.1 is called the {\it Z$_3$-graded Leibniz rule}.

For any differential calculus over the Z$_3$-graded algebra the following identity holds:
\begin{eqnarray} \label{11}
{\sf d}(a \, {\sf d}b) = {\sf d}a\wedge {\sf d}b + q^{p(a)} \, a \, {\sf d}^2b.
\end{eqnarray}
We recall here that $q=j=e^{2\pi i/3}$.

\begin{defin}
A differential calculus $\Gamma^\wedge$ over the Z$_3$-graded quantum superspace ${\mathbb R}_q(2|1)$ with left coaction $\varphi:{\mathbb R}_q(2|1) \longrightarrow {\cal A}\otimes{\mathbb R}_q(2|1)$ is called left covariant with respect to ${\cal A}$ if there exists an algebra homomorphism $\Delta_L:\Gamma^\wedge\longrightarrow {\cal A}\otimes\Gamma^\wedge$ which is a left coaction of ${\cal A}$ on $\Gamma^\wedge$ such that

(1) $\Delta_L(a\rho b) = \varphi(a)\Delta_L(\rho)\varphi(b)$ for all $a,b \in {\mathbb R}_q(2|1)$ and $\rho\in\Gamma^\wedge$,

(2) $\Delta_L({\sf d}a) = (\tau\otimes {\sf d})\varphi(a)$ for all $a\in{\mathbb R}_q(2|1)$,

(3) $\Delta_L({\sf d}\rho) = (\tau\otimes {\sf d})\Delta_L(\rho)$ for $\rho\in\Gamma^\wedge$

\noindent
where $\tau:\Gamma^\wedge\longrightarrow\Gamma^\wedge$ is the linear map of degree zero which gives $\tau(\rho)=q^{p(\rho)}\rho$ for all $\rho\in\Gamma^\wedge$.
\end{defin}
For the Hopf algebra ${\cal A}$ presented in Section 2, we will take $\varphi=\Delta$. From the Definition 3.2 one deduces the following properties \cite{Woro}:
\begin{eqnarray} \label{12}
(\Delta\otimes\mbox{id})\circ\Delta_L = (\mbox{id}\otimes\Delta_L)\circ\Delta_L, \quad
(\epsilon\otimes\mbox{id})\Delta_L({\sf d}u) = {\sf d}u.
\end{eqnarray}

It is well known that in classical differential calculus, functions commute with differentials. From an algebraic point of view, the space of 1-forms is a free finite bimodule over the algebra of smooth functions generated by the first order differentials and the commutativity shows how its left and right structure are related to each other.

Note that the linear operator {\sf d} applied to $x$ produces a 1-form whose Z$_3$-grade is 1, by definition. Similarly, application of {\sf d} to $y$ and $\theta$ produces  1-forms whose Z$_3$-grade are 2 and 0, respectively. We shall denote the obtained quantities by ${\sf d}x$, ${\sf d}y$ and ${\sf d}\theta$. When the linear operator {\sf d} is applied to the 1-forms, it will produce new entities which we shall call 2-forms of grade 2, 0, 1, respectively, denoted by ${\sf d}^2x$, and so on. Finally, we demand that ${\sf d}^3=0$.

\subsection{Structure of the commutation relations} \label{sec3.2}

To define a differential geometry on the algebra ${\cal A}$ one has to construct, first of all, an analog of a differential calculus on the algebra. To this aim let us associate generators $x$, $y$ and $\theta$ and their differentials ${\sf d}x$, ${\sf d}y$ and ${\sf d}\theta$ which are considered as generators of a space $\Gamma^1({\cal A})\doteq \Gamma^1$ of 1-forms. We allow a multiplication of the differentials by the elements of ${\cal A}$ from the left and from the right and by the definition of the multiplications the resulting 1-form belongs to $\Gamma^1$ again. This means that $\Gamma^1$ is an ${\cal A}$-bimodule.

By the condition (2) of the Definition 3.2, we have
\begin{eqnarray} \label{13}
\Delta_L(x) &=& x \otimes {\sf d}x, \nonumber\\
\Delta_L(y) &=& x \otimes {\sf d}y, \nonumber\\
\Delta_L(\theta) &=& q^2 \, \theta \otimes {\sf d}x + x^{-1} \otimes {\sf d}\theta.
\end{eqnarray}

Now, in order to obtain the Z$_3$-graded commutation relations of the elements of the algebra ${\cal A}$ with their differentials, we shall use the approach of \cite{WZ} and \cite{Celik2}.

\begin{thm} 
The Z$_3$-graded commutation relations satisfied by the generators of ${\cal A}$ and their first order differentials are in the form
\begin{eqnarray} \label{14}
x \, {\sf d}x=q \, {\sf d}x \, x, & x \, {\sf d}y = q^2 \, {\sf d}y \, x, & x \, {\sf d}\theta=q \, {\sf d}\theta \, x+(q-1) \,{\sf d}x \, \theta, \nonumber \\
y \, {\sf d}y={\sf d}y \, y, & y \, {\sf d}x = q \, {\sf d}x \, y + (1-q^2) \,{\sf d}y \, x, & y \, {\sf d}\theta = q \, {\sf d}\theta \, y + (1-q^2) \,{\sf d}y \, \theta, \nonumber\\
\theta \, {\sf d}x = q \, {\sf d}x \, \theta, & \theta \, {\sf d}y = {\sf d}y \, \theta, & \theta \, {\sf d}\theta = q^2 \, {\sf d}\theta \, \theta.
\end{eqnarray}
\end{thm}

\begin{pf}
To obtain these relations, we assume that the possible commutation relations of the generators with their differentials are of the following form
\begin{eqnarray} \label{15}
x \,{\sf d}x &=& Q_1 \, {\sf d}x \,  x, \quad x \, {\sf d}y = A_{11} \, {\sf d}y \, x + A_{12} \, {\sf d}x \,  y, \nonumber \\
y \, {\sf d}y &=& Q_2 \, {\sf d}y \,  y, \quad y \, {\sf d}x = A_{21} \, {\sf d}x \,  y + A_{22} \, {\sf d}y \, x, \nonumber \\
\theta \, {\sf d}\theta &=& Q_3 \, {\sf d}\theta\,\theta, \quad x \, {\sf d}\theta = B_{11} \, {\sf d}\theta \, x + B_{12} \, {\sf d}x \, \theta, \nonumber \\
\theta \, {\sf d}x &=& B_{21} \, {\sf d}x \, \theta + B_{22} \, {\sf d}\theta \, x, \quad y \, {\sf d}\theta = C_{11} \, {\sf d}\theta \, y + C_{12} \, {\sf d}y \, \theta, \nonumber \\
\theta \, {\sf d}y &=& C_{21} \, {\sf d}y \, \theta + C_{22} \, {\sf d}\theta \, y
\end{eqnarray}
where the coefficients $Q_i$, $A_{ij}$, $B_{ij}$, $C_{ij}$ are possible related to $q$. One knows, from \cite{WZ}, that the application of the exterior differential {\sf d} to (\ref{15}) does not give commutation relations between the differentials ${\sf d}x$, ${\sf d}y$ and ${\sf d}\theta$. Therefore, to find them we first apply the operator $\Delta_L$ to the relations (\ref{15}). If we consider the first four relations in (\ref{15}), and we apply them to the operator $\Delta_L$, we find $A_{12}=0$, $A_{21}=q$ and $Q_2=1$. Next, to obtain the remaining coefficients, we use the consistency of calculus. The consistency conditions give rise to $(A_{11}-qQ_1)A_{22}=0$ and $A_{22}=qA_{11}-q^2$. Thus, we have four solutions depending on the choice of $Q_1$. Let's take as an appropriate choice $Q_1=q$ and $A_{22}\ne0$. Then we get the first, second, fourth and fifth relations in (\ref{14}). When similar operations are performed, the other relations in (\ref{14}) can be also obtained.
\end{pf}

\begin{thm} 
Commutation relations between the first order differentials of the generators of ${\cal A}$ are
\begin{eqnarray} \label{16}
{\sf d}x\wedge{\sf d}y = q^2 \,{\sf d}y\wedge{\sf d}x, \quad {\sf d}x\wedge{\sf d}\theta = {\sf d}\theta\wedge{\sf d}x, \quad
{\sf d}y\wedge{\sf d}\theta = {\sf d}\theta\wedge{\sf d}y.
\end{eqnarray}
\end{thm}

\begin{pf}
To obtain above relations, we assume that the possible commutation relations of the differentials are of the following form
\begin{eqnarray} \label{17}
{\sf d}x \wedge {\sf d}y = F_1 \, {\sf d}y \wedge {\sf d}x, \quad {\sf d}x \wedge {\sf d}\theta = F_2 \, {\sf d}\theta \wedge {\sf d}x, \quad
{\sf d}y \wedge {\sf d}\theta = F_3 \, {\sf d}\theta \wedge {\sf d}y,
\end{eqnarray}
where the coefficients $F_1, F_2, F_3$ are possible related to $q$. The actions of the operator $\Delta_L$ on the relations (\ref{17}) give rise $F_2=1$ and $F_3=q^2F_1^{-1}$. There are no restrictions on the parameter $F_1$, but it will be appropriate to choose it as $q^2$.
\end{pf}

\noindent{\bf Note 3.} Using the first and fourth relations in (\ref{14}), it can easily be shown that
\begin{equation*} \label{}
{\sf d}({\sf d}x\wedge{\sf d}x\wedge{\sf d}x) = 0 = {\sf d}({\sf d}y\wedge{\sf d}y\wedge{\sf d}y).
\end{equation*}
For the first one, we apply twice the operator {\sf d} to the first relation in (\ref{14}). The result will give the relation ${\sf d}x\wedge{\sf d}^2x = q \,{\sf d}^2x\wedge{\sf d}x$ with the assumption ${\sf d}^3=0$. Now using this relation, we get the first equation above. So, we can choose
$$\label{16'}
{\sf d}x\wedge{\sf d}x\wedge{\sf d}x = 0 \quad \mbox{and} \quad {\sf d}y\wedge{\sf d}y\wedge{\sf d}y = 0. \eqno(16')$$

If we introduce the differentials of the coordinate functions as
\begin{equation*}
\xi = {\sf d}x, \quad \eta = {\sf d}y, \quad z = {\sf d}\theta
\end{equation*}
the relations (\ref{2}) coincide with (\ref{16}) and (16$'$).

Since ${\sf d}^3=0$ (and ${\sf d}^2\neq 0$) in the Z$_3$-graded space, in order to construct a self-consistent theory of differential forms it is necessary to add to the first order differentials of coordinates ${\sf d}x,{\sf d}y,{\sf d}\theta$ a set of second order differentials ${\sf d}^2x, {\sf d}^2y, {\sf d}^2\theta$. Appearance of higher order differentials is a peculiar property of a proposed generalization of differential forms. Now we will get the commutation relations of the generators $x$, $y$ and $\theta$ with their second order differentials ${\sf d}^2x$, ${\sf d}^2y$ and ${\sf d}^2\theta$ which are considered as generators of a space $\Gamma^2({\cal A})\doteq \Gamma^2$ of 2-forms.

\begin{thm} 
The commutation relations of the generators of the algebra ${\cal A}$ and their second order differentials are in the form
\begin{eqnarray} \label{18}
x \, {\sf d}^2x  &=& q \, {\sf d}^2x \, x + (q^2-1) \, {\sf d}x\wedge{\sf d}x, \quad x \, {\sf d}^{2}y  = q^2 \, {\sf d}^2y \, x + (q^2-1) \, {\sf d}x\wedge{\sf d}y, \nonumber \\
x \, {\sf d}^2\theta &=& q \, {\sf d}^2\theta \, x + (q-1) \, {\sf d}^2x \, \theta + (q^2-1) \, {\sf d}x\wedge{\sf d}\theta, \nonumber \\
y \, {\sf d}^2x  &=& {\sf d}^2x \, y + (q^2-q) \, {\sf d}^2y \, x + (q-q^2) \, {\sf d}y\wedge{\sf d}x, \\
y \, {\sf d}^2y &=& q^2 \, {\sf d}^2y \, y + (q-q^2) \, {\sf d}y\wedge{\sf d}y, \nonumber \\
y \, {\sf d}^2\theta &=& {\sf d}^2\theta \, y + (q^2-q) \, {\sf d}^2y \, \theta + (q-q^2) \, {\sf d}y\wedge{\sf d}\theta, \nonumber \\
\theta \, {\sf d}^2x &=& q^2 \, {\sf d}^2x \, \theta + (1-q) \, {\sf d}\theta\wedge{\sf d}x, \quad \theta \, {\sf d}^2y = q \, {\sf d}^2y \, \theta + (1-q) \, {\sf d}\theta\wedge{\sf d}y, \nonumber \\
\theta \, {\sf d}^2\theta &=& {\sf d}^2\theta \, \theta + (1-q) \, {\sf d}\theta\wedge{\sf d}\theta. \nonumber
\end{eqnarray}
\end{thm}

\begin{pf}
To obtain the relations (\ref{18}) we will apply the exterior differential {\sf d} from left to the relations (\ref{14}). If we differentiate the relations in (\ref{14}) with respect to the Z$_3$-graded Leibniz rule we get the relations (\ref{18}) provided that $F_1=q^2$.
\end{pf}
It is easily seen that the linear map $\Delta_L$ leaves invariant the relations (\ref{18}) in the sense that
\begin{eqnarray*}
\Delta_L(u \, {\sf d}^2v) = \Delta(u)(\tau\otimes{\sf d}^2)\Delta(v).
\end{eqnarray*}

The relations (\ref{18}) are not homogeneous in the sense that the commutation relations between the generators and their second order differentials include first order differentials as well.

We now proceed in the construction of the differential calculus by introducing the commutation relations of the first order differentials with second order differentials.

\begin{thm} 
The commutation relations between the first order differentials and the second order differentials are as follows
\begin{eqnarray} \label{19}
{\sf d}x \wedge {\sf d}^2x &=&  q \, {\sf d}^2x \wedge {\sf d}x, \quad {\sf d}x \wedge {\sf d}^2y = q \, {\sf d}^2y \wedge {\sf d}x + (q-q^2) \, {\sf d}^2x \wedge {\sf d}y, \nonumber \\
{\sf d}x \wedge {\sf d}^2\theta &=& q \, {\sf d}^2\theta \wedge {\sf d}x, \quad
{\sf d}y \wedge {\sf d}^2x = q \, {\sf d}^2x \wedge {\sf d}y, \quad {\sf d}y \wedge {\sf d}^2y = {\sf d}^2y \wedge {\sf d}y,\nonumber \\
{\sf d}y \wedge {\sf d}^2\theta &=& {\sf d}^2\theta \wedge {\sf d}y, \quad {\sf d}\theta \wedge {\sf d}^2x = {\sf d}^2x \wedge {\sf d}\theta + (q^2-1) \, {\sf d}^2\theta \wedge {\sf d}x, \\
{\sf d}\theta \wedge {\sf d}^2\theta &=& q^2 \, {\sf d}^2\theta \wedge {\sf d}\theta, \quad
{\sf d}\theta \wedge {\sf d}^2y = {\sf d}^2y \wedge {\sf d}\theta + (q^2-1) \, {\sf d}^2\theta \wedge {\sf d}y. \nonumber
\end{eqnarray}
\end{thm}

\begin{pf}
If we apply the exterior differential {\sf d} to the relations (\ref{18}) using assumption that ${\sf d}^3=0$ then we get the relations (\ref{19}).
\end{pf}

\noindent{\bf Note 4.} In Note 3, we chose ${\sf d}x\wedge{\sf d}x\wedge{\sf d}x=0$ and ${\sf d}y\wedge{\sf d}y\wedge{\sf d}y=0$. However, the quantities ${\sf d}x\wedge{\sf d}x$ and ${\sf d}y\wedge{\sf d}y$ do not have to be zero. Indeed,
\begin{equation*}
{\sf d}({\sf d}x\wedge{\sf d}x)=-{\sf d}x\wedge{\sf d}^2x \quad \mbox{and} \quad {\sf d}({\sf d}y\wedge{\sf d}y)=-q \, {\sf d}y\wedge{\sf d}^2y.
\end{equation*}

The relations in the following proposition can be obtained with applying the exterior differential {\sf d} to the relations (\ref{19}).

\begin{thm} 
The commutation relations between the second order differentials are as follows
\begin{eqnarray} \label{20}
{\sf d}^2x \wedge {\sf d}^2y &=& {\sf d}^2y \wedge {\sf d}^2x, \nonumber \\
{\sf d}^2y \wedge {\sf d}^2\theta &=& q \, {\sf d}^2\theta \wedge {\sf d}^2y, \nonumber \\
{\sf d}^2x \wedge {\sf d}^2\theta &=& q^2 \, {\sf d}^2\theta \wedge {\sf d}^2x.
\end{eqnarray}
\end{thm}

Consequently, we set up a differential schema with the relations satisfied by the elements of the set $\{x,y,\theta,{\sf d}x,{\sf d}y,{\sf d}\theta,{\sf d}^2x,{\sf d}^2y,{\sf d}^2\theta\}$.

\subsection{Cartan-Maurer one forms on $\cal A$ }\label{sec3.3}

We now proceed in the construction of the differential calculus by introducing the space of 1-forms. We shall define three 1-forms using the generators of ${\cal A}$ and investigate their relations with the coordinates and themselves. Let us begin the following definition given by Woronowicz \cite{Woro}.

\begin{defin}
A left-covariant bimodule over the Hopf algebra ${\cal A}$ is an ${\cal A}$-bimodule $\Gamma$ which is a left comodule of ${\cal A}$ with left coaction $\Delta_L:\Gamma\longrightarrow {\cal A}\otimes\Gamma$, such that
\begin{eqnarray*}
\Delta_L(a\rho b) &=& \Delta(a)\Delta_L(\rho)\Delta(b)
\end{eqnarray*}
for all $a,b \in {\cal A}$ and $\rho\in\Gamma$,
\end{defin}
For a left-covariant bimodule $\Gamma^\wedge$, an element of the vector space
\begin{eqnarray*}
\Omega_L = \{\rho\in \Gamma^\wedge: \Delta_L(\rho) = 1\otimes \rho\}
\end{eqnarray*}
is called {\it left-invariant}.

Since $\Gamma^\wedge$ is a left-covariant differential calculus over ${\cal A}$, one can define a linear mapping (the Cartan-Maurer 1-form) $w_a:{\cal A}\longrightarrow\Omega_L$ by setting \cite{Woro}
\begin{eqnarray} \label{21}
w_{\Gamma^\wedge}(a) &=& (\kappa\otimes{\sf d})\Delta(a).
\end{eqnarray}
If we omit the subscript $\Gamma^\wedge$ and write simply $w_a$, then we have
\begin{eqnarray} \label{22}
w_x = x^{-1} \, {\sf d}x, \quad w_y = x^{-1} \, {\sf d}y, \quad w_\theta = x \, {\sf d}\theta - \theta \, {\sf d}x.
\end{eqnarray}
Here, the 1-forms $w_x$, $w_y$, $w_\theta$ are of degree 1, 2, 0, respectively. It follows immediately, from (\ref{4}) and (\ref{13}), that $\Delta_L(w_a)=1\otimes w_a$ for all $a\in {\cal A}$. So, we get $w({\cal A})=\Omega_L$.

The commutation relations between the generators of ${\cal A}$ and their differentials can be expressed in terms of the Cartan-Maurer 1-forms:

\begin{thm} 
The commutation relations of above 1-forms with the generators of ${\cal A}$ are as follows
\begin{eqnarray} \label{23}
x \, w_x &=& q \, w_x \, x, \quad x \, w_y = q^2 \, w_y \, x, \quad x \, w_\theta = q \, w_\theta \, x, \nonumber \\
y \, w_x &=& q^2 \, w_x \, y + (q-1) \, w_y \, x, \quad y \, w_y = q \, w_y \, y, \quad y \, w_\theta = w_\theta \, y, \nonumber \\
\theta \, w_x &=& q^2 \, w_x \, \theta, \quad \theta \, w_y = q \, w_y \, \theta, \quad \theta \, w_\theta = q \, w_\theta \, \theta.
\end{eqnarray}
\end{thm}

Using (\ref{14}) and (\ref{16}) we now find the commutation rules of the generators of $\Omega_L$ as follows

\begin{thm} 
The commutation rules of the elements of $\Omega_L$ are as follows
\begin{eqnarray} \label{24}
w_x \wedge w_y &=& q \, w_y \wedge w_x, \nonumber\\
w_y \wedge w_\theta  &=& w_\theta \wedge w_y,  \nonumber\\
w_\theta \wedge w_x &=& q^2 \, w_x \wedge w_\theta, \nonumber\\
w_x\wedge w_x\wedge w_x &=& 0 = w_y\wedge w_y\wedge w_y.
\end{eqnarray}
\end{thm}

Later (in the proof of Proposition 4.2) we need the explicit formulae for ${\sf d}^2w_x$, ${\sf d}^2w_y$ and ${\sf d}^2w_\theta$. From (\ref{22}), one can write
\begin{equation*} \label{}
{\sf d}w_x = x^{-1} \, {\sf d}^2x - w_x\wedge w_x, \quad {\sf d}w_y = x^{-1} \, {\sf d}^2y - w_x\wedge w_y, \quad {\sf d}w_\theta = x \, {\sf d}^2\theta - q^2 \,\theta \, {\sf d}^2x.
\end{equation*}
Remembering that ${\sf d}^2\ne0$ and ${\sf d}^3=0$, we get the following 3-forms
\begin{eqnarray*} \label{}
{\sf d}^2w_x &=& 0, \\
{\sf d}^2w_y &=& q^2 \, {\sf d}w_y\wedge w_x - q \, {\sf d}w_x\wedge w_y, \\
{\sf d}^2w_\theta &=& q \, {\sf d}w_\theta\wedge w_x - {\sf d}w_x\wedge w_\theta - q^2 \, w_x\wedge w_x\wedge w_\theta.
\end{eqnarray*}

\subsection{Relations of partial derivatives}\label{sec3.4}
We now proceed in the construction of the differential calculus by introducing the partial derivatives of the generators of the algebra ${\cal A}$ and we will obtain the commutation relations between the coordinates, first and second order differentials and their partial derivatives.

Since $\Gamma^\wedge$ is a left covariant differential calculus, for any $a\in{\cal A}$ there are uniquely determined elements $\partial_i(a)\in{\cal A}$ such that
\begin{eqnarray*}
{\sf d}a &=& {\sf d}x \, \partial_x(a) + {\sf d}y \, \partial_y(a) + {\sf d}\theta \, \partial_\theta(a).
\end{eqnarray*}
The mappings $\partial_i:{\cal A}\longrightarrow{\cal A}$ are called the {\it partial derivatives} of the calculus $\Gamma^\wedge$. The following proposition gives the relations between the elements of ${\cal A}$ and their partial derivatives.

\begin{thm} 
The relations of the coordinates with their partial derivatives are as follows
\begin{eqnarray} \label{25}
\partial_x \, x &=& 1+q \, x \, \partial_x + (q-1) \, \theta \, \partial_\theta, \quad \partial_x \, y = q^2 \, y \, \partial_x, \quad \partial_x \, \theta=\theta \, \partial_x, \nonumber \\
\partial_y \, x &=& q^2 \, x \, \partial_y, \quad \partial_y \, \theta = q^2 \, \theta \, \partial_y, \\
\partial_y \, y &=& 1 + q \, y \, \partial_y + (q-1) \, (x \, \partial_x + \theta \, \partial_\theta), \nonumber \\
\partial_\theta \, x &=&  q \, x \, \partial_\theta, \quad \partial_\theta \, y = q^2 \, y \, \partial_\theta, \quad
\partial_\theta \, \theta = 1 + q \, \theta \, \partial_\theta. \nonumber
\end{eqnarray}
\end{thm}

\begin{pf}
We know that, the exterior differential {\sf d} can be written in terms of the differentials and partial derivatives as
\begin{eqnarray} \label{26}
{\sf d}f &=& ({\sf d}x \, \partial_x + {\sf d}y \, \partial_y + {\sf d}\theta \, \partial_\theta)f
\end{eqnarray}
where $f$ is a differentiable function. For consistency, the grades of the partial derivatives $\partial_x$, $\partial_y$ and $\partial_\theta$ should be 0, 2, 1, respectively. Now, if we replace $f$ with $xf$ in the left hand side of the equality in (\ref{26}), we get
\begin{eqnarray*}
{\sf d}(x f)
&=& {\sf d}x \, f+ x \, {\sf d}f \nonumber\\
&=& \left[{\sf d}x \, (1+q \, x \, \partial_x + (q-1) \, \theta \, \partial_\theta) + {\sf d}y \, (q^2 \, x \, \partial_y) + {\sf d}\theta \, (q \, x \, \partial_\theta)\right]f. \nonumber
\end{eqnarray*}
On the other hand, the right hand side of the equality in (\ref{26}) has the form
\begin{eqnarray*}
{\sf d}(xf) &=& [{\sf d}x \, (\partial_x \, x) + {\sf d}y \, (\partial_y \, x) + {\sf d}\theta \, (\partial_\theta \, x)]f.
\end{eqnarray*}
Now, by comparing the right hand sides of these two equalities according to the differentials we obtain some of the relations in (\ref{25}). Other relations can be found similarly.
\end{pf}

\begin{thm} 
The relations between partial derivatives are in the form
\begin{eqnarray} \label{27}
\partial_x \, \partial_y &=& q^2 \, \partial_y \, \partial_x, \nonumber \\
\partial_y \, \partial_\theta &=& q \, \partial_\theta \, \partial_y, \nonumber \\
\partial_\theta \, \partial_x  &=& \partial_x \, \partial_\theta \, , \qquad \partial_\theta^3 = 0.
\end{eqnarray}
\end{thm}

\begin{pf}
To obtain the relations in (\ref{27}) one uses ${\sf d}^3=0$ with the relations given in (\ref{18}) and (\ref{20}).
\end{pf}

\noindent{\bf Note 5.} It is easily seen that the cube of the partial derivatives $\partial_x$ and $\partial_y$ are both central.
We said above that, the partial derivative $\partial_y$ is of degree 2. Therefore, in terms of consistency, we can assume that $\partial_y^3=0$.

\begin{defin} 
The Z$_3$-graded quantum Weyl superalgebra ${\cal W}_q(2|1)$ is the unital algebra generated by $x$, $y$, $\theta$ and $\partial_x$, $\partial_y$, $\partial_\theta$ which satisfy the relations (\ref{1}), (\ref{25}) and (\ref{27}).
\end{defin}

To complete the schema, we need the relations that partial derivatives and first and second order differentials satisfied.

\begin{thm} 
The relations partial derivatives with first and second order differentials are as follows
\begin{eqnarray} \label{28}
\partial_x \, {\sf d}x  &=& q^2 \, {\sf d}x \, \partial_x + (q^2-1) \, {\sf d}y \, \partial_y, \quad \partial_x \, {\sf d}y = q \, {\sf d}y \, \partial_x, \quad \partial_x \, {\sf d}\theta = q^2 \, {\sf d}\theta \, \partial_x,\nonumber \\
\partial_y \, {\sf d}x &=& q^2 \, {\sf d}x \, \partial_y, \quad \partial_y \, {\sf d}y = {\sf d}y \, \partial_y, \quad
\partial_y \, {\sf d}\theta = q^2 \, {\sf d}\theta \, \partial_y, \\
\partial_\theta \,{\sf d}x &=& q^2 \,{\sf d}x \,\partial_\theta, \,\,
\partial_\theta \,{\sf d}y = {\sf d}y \,\partial_\theta, \quad \partial_\theta \, {\sf d}\theta = q \,{\sf d}\theta \, \partial_\theta + (q-q^2)({\sf d}x \, \partial_x + {\sf d}y \, \partial_y)\nonumber
\end{eqnarray}
and
\begin{eqnarray} \label{29}
\partial_x \, {\sf d}^2x &=&q^2 \, {\sf d}^2x \, \partial_x+(q^2-1) \, {\sf d}^2y \, \partial_y, \quad \partial_x \, {\sf d}^2y=q \, {\sf d}^2y \,  \partial_x,\\
\partial_x \, {\sf d}^2\theta &=& q^2 \, {\sf d}^2\theta \, \partial_x, \,\,
\partial_y \,  {\sf d}^2x = {\sf d}^2x \, \partial_y, \quad \partial_y \, {\sf d}^2y = q \, {\sf d}^2y \, \partial_y, \quad
\partial_y \, {\sf d}^2\theta = {\sf d}^{2}\theta \, \partial_y, \nonumber \\
\partial_\theta \, {\sf d}^2x &=& q^2{\sf d}^2x \, \partial_\theta, \,\,\, \partial_\theta \, {\sf d}^2y = q^2{\sf d}^2y \, \partial_\theta, \quad
\partial_\theta \, {\sf d}^2\theta = q^2{\sf d}^2\theta \, \partial_\theta +(q^2-q){\sf d}^2x \, \partial_x.\nonumber
\end{eqnarray}
\end{thm}

\begin{pf}
To obtain the relations (\ref{28}), we assume that the relations between partial derivatives and first order differentials as follows
\begin{eqnarray*}
\partial_x \, {\sf d}x &=& A_1 \, {\sf d}x \, \partial_x + A_2 \, {\sf d}y \, \partial_y + A_3 \, {\sf d}\theta \, \partial_\theta,  \\
\partial_x \, {\sf d}y &=& A_4 \, {\sf d}y \, \partial_x + A_5 \, {\sf d}x \, \partial_y, \\
\partial_x \, {\sf d}\theta &=& A_6 \, {\sf d}\theta \, \partial_x + A_7 \, {\sf d}x \, \partial_\theta, \\
&& \vdots
\end{eqnarray*}
and so on. The coefficients $A_i$'s will be determined in terms of the deformation parameter $q$. To find them we apply the operator  $\partial_x$ to the relations (\ref{14}). Using the relation
\begin{eqnarray*}
\partial_i (X^j {\sf d}X^k) = \delta^i{}_j \delta^k{}_l {\sf d}X^k
\end{eqnarray*}
where $\partial_1 = \partial_x$, $\partial_2 = \partial_y$, $\partial_3 = \partial_\theta$, $X^1 = x$, $X^2=y$ and $X^3 = \theta$
after from some calculations we find $A_1=q^2 \, $, $A_2=q^2-1$, $A_3=0$, $A_4=q$, $A_5=0$, $A_6=q^2$ and $A_7=0$. Other relations in (\ref{28}) can be determined in a similar manner. The relations in (\ref{29}) can also be obtained using the same idea.
\end{pf}

\noindent{\bf Note 6.} Using the relations (\ref{28}), it can be easily seen that the relations
\begin{eqnarray*}
\partial_i \, {\sf d} = q^{2-p(\partial_i)} \, {\sf d} \, \partial_i \quad \mbox{and} \quad \partial_i \, {\sf d}^2 = q^{1+p(\partial_i)} \, {\sf d}^2 \, \partial_i
\end{eqnarray*}
are consistent with ${\sf d}^3=0$.

\section{Z$_3$-graded quantum Lie superalgebra}\label{sec4}

The commutation relations of Cartan-Maurer forms allow us to construct the algebra of the generators. In order to obtain the Z$_3$-graded quantum Lie superalgebra of the algebra generators we first write the Cartan-Maurer forms as
\begin{eqnarray} \label{30}
{\sf d}x = x \, w_x, \quad {\sf d}y = x \, w_y, \quad {\sf d}\theta = x^{-1} \, w_\theta + q^2 \, \theta \, w_x.
\end{eqnarray}
Then the differential ${\sf d}$ can be expressed in the form
\begin{equation} \label{31}
{\sf d}f = (w_x \,  T_x + w_y \,  T_y + w_\theta \,  T_\theta)f.
\end{equation}
Here $T_x$, $T_y$ and $T_\theta$ are the (quantum) Lie superalgebra generators which are of degree 0,2,1, respectively. We now shall give the commutation relations of these generators.

\begin{lemma} 
Commutation relations between 1-forms and their differentials are as follows
\begin{eqnarray*}
w_x\wedge {\sf d}w_x &=& q \, {\sf d}w_x \wedge w_x, \noindent\\
w_x\wedge {\sf d}w_y &=& {\sf d}w_y \wedge w_x + (q-q^2) ({\sf d}w_x + w_x \wedge w_x)\wedge w_y, \noindent\\
w_x\wedge {\sf d}w_\theta &=& q^2 \, {\sf d}w_\theta \wedge w_x + (q-q^2) ({\sf d}w_x - w_x \wedge w_x)\wedge w_\theta, \noindent\\
w_y\wedge {\sf d}w_x &=& q^2 \, {\sf d}w_x \wedge w_y + (q^2-1) \, w_x \wedge w_x \wedge w_y, \noindent\\
w_y\wedge {\sf d}w_y &=& {\sf d}w_y \wedge w_y + (1-q) \, w_x \wedge w_y \wedge w_y, \noindent\\
w_y\wedge {\sf d}w_\theta &=& {\sf d}w_\theta \wedge w_y + (1-q^2) \, w_x \wedge w_y \wedge w_\theta, \noindent\\
w_\theta\wedge {\sf d}w_x &=& q \, {\sf d}w_x \wedge w_\theta + (q^2-1) \,({\sf d}w_\theta - w_x \wedge w_\theta)\wedge w_x, \noindent\\
w_\theta\wedge {\sf d}w_y &=& {\sf d}w_y \wedge w_\theta + (q^2-q) \,(q \,{\sf d}w_\theta - w_x \wedge w_\theta)\wedge w_y, \noindent\\
w_\theta\wedge {\sf d}w_\theta &=& q^2 \, {\sf d}w_\theta \wedge w_\theta.
\end{eqnarray*}
\end{lemma}

The proof of this lemma is easy and straightforward, but you will have to do some tedious process.

\begin{thm} 
The operators $T_x$, $T_y$ and $T_\theta$ constitute a Z$_3$-graded quantum Lie superalgebra with the following Lie brackets
\begin{eqnarray} \label{32}
\left[T_x, T_y\right]_q &=& q \,T_y, \nonumber \\
\left[T_x, T_\theta\right]_{q} &=& q \,T_\theta, \quad T_\theta^3 = 0, \nonumber \\
\left[T_y, T_\theta\right]_q &=& 0, \qquad T_y^3 = 0.
\end{eqnarray}
\end{thm}

\begin{pf}
The assumption ${\sf d}^3=0$ and long calculations with use of Lemma 4.1 give rise the commutation relations in (\ref{32}). However, it is easy to show that the operators $T_x$, $T_y$ and $T_\theta$ satisfy the Jacobi identity.
\end{pf}

The commutation relations (\ref{32}) of the Lie superalgebra generators should be consistent with monomials of the coordinates of the Z$_3$-graded quantum  superspace. For this, we evaluate the commutation relations between the generators of superalgebra and the coordinates. These relations can be extracted from the Z$_3$-graded Leibniz rule:

\begin{lemma} 
The commutation relations of the generators with the coordinates are as follows
\begin{eqnarray} \label{33}
T_x \, x &=& q \, x + q \, x \, T_x, \quad T_x \, y = y \, T_x, \quad T_x \, \theta = q \, \theta + q \, \theta \, T_x, \nonumber\\
T_y \, x &=& q^2 \, x \, T_y, \quad T_y \, y = q^2 \, x + q^2 \, y \, T_y + (q^2-q) \, x \, T_x,  \quad T_y \, \theta = \theta \, T_y, \nonumber \\
T_\theta \, x &=& q \, x \, T_\theta, \quad T_\theta \, y = q^2 \, y \, T_\theta, \quad T_\theta \, \theta = q^2 \, x^{-1} + \theta \, T_\theta.
\end{eqnarray}
\end{lemma}

The following corollary can be proven by using the Mathematical Induction Principle.
\begin{corollary} 
For $k \in {\mathbb N}_0$, one has
\begin{eqnarray*}
T_x \, x^k = \frac{1-q^k}{1-q^2} \, x^k + q^k x^k \, T_x, \quad T_y \, x^k = q^{2k} x^k \, T_y, \quad T_\theta \, x^k = q^k x^k \, T_\theta.
\end{eqnarray*}
\end{corollary}

The proof of the following proposition is based on Corollary 4.4.

\begin{thm} 
The coproduct structure of the Lie algebra generators is given by
\begin{eqnarray} \label{34}
\Delta(T_x) & = & T_x \otimes {\bf 1} + \left[{\bf 1} + (1-q^2) T_x\right] \otimes T_x, \nonumber \\
\Delta(T_y) & = & T_y \otimes \left[{\bf 1} + (1-q^2) T_x\right]^2 + {\bf 1} \otimes T_y, \nonumber \\
\Delta(T_\theta) & = & T_\theta \otimes \left[{\bf 1} + (1-q^2) T_x \right] + \left[{\bf 1} + (1-q^2) T_x \right]^2\otimes T_\theta.
\end{eqnarray}
\end{thm}

The action of the counit and the coinverse on the vector fields can be obtained from the identities given in (\ref{8}) and (\ref{9}).

We know, from subsection 3.4, that the exterior differential {\sf d} can be expressed in the form (\ref{26}), that is,
\begin{equation*}
{\sf d}f = ({\sf d}x \, \partial_x + {\sf d}y \, \partial_y + {\sf d}\theta \, \partial_\theta)f.
\end{equation*}
Considering (\ref{26}) together (\ref{31}) and using (\ref{22}) one has
\begin{equation} \label{35}
T_x \equiv q \, (x \, \partial_x + \theta \, \partial_\theta), \quad  T_y \equiv q^2 \, x \, \partial_y, \quad T_\theta \equiv q^2 \, x^{-1} \, \partial_\theta. \end{equation}
Using the relations (\ref{25}) and (\ref{27}) one can check that the relation of the generators in (\ref{35}) coincide with (\ref{32}). It can also be verified that, the action of the generators in (\ref{35}) on the coordinates coincide with (\ref{33}).

\section{The dual of the Hopf algebra ${\cal A}$}\label{sec5}
Let us begin with the Hopf algebra ${\cal A}$. Then its dual ${\cal U}\doteq {\cal A}'$ is a Hopf algebra as well.
Using the coproduct $\Delta$ in ${\cal A}$, one defines a product in ${\cal U}$ and using the product in the Hopf algebra ${\cal A}$, one defines a coproduct
in ${\cal U}$.

In this section, in order to obtain the dual of the Hopf algebra ${\cal A}$ defined in section 2, we have applied
the approach which Sudbery invented for $A_q(2)$ \cite{SD} to ${\cal A}=O({\mathbb R}_q(2|1))$. Let us begin the definition of the duality given by Abe \cite{Abe}.

\begin{defin} 
A dual pairing of two bialgebras ${\cal U}$ and ${\cal A}$ is a bilinear mapping $<,>:{\cal U}\times{\cal A}\longrightarrow {\mathbb C}$ such that
\begin{eqnarray} \label{36}\label{37}
<u, ab> &=& <\Delta_{\cal U}(u), a \otimes b>, \\
<uv, a> &=& <u \otimes v, \Delta_{\cal A}(a)>, \\
<u, 1_{\cal A}> &=& \epsilon_{\cal U}(u), \quad <1_{\cal U}, a> = \epsilon_{\cal A}(a) \nonumber
\end{eqnarray}
for all $u, v \in {\cal U}$ and $a, b \in {\cal A}$. We say that the pairing is non-degenerate if
$<u,a>=0$ for all $a \in {\cal A}$ implies that $u = 0$, and $<u,a>=0$ for all $u \in {\cal U}$ implies that $a=0$.
\end{defin}

Two Hopf algebras ${\cal U}$ and ${\cal A}$ are said to be in duality if they are in duality as bialgebras and if
\begin{eqnarray} \label{38}
<\kappa_{\cal U}(u), a> &=& <u, \kappa_{\cal A}(a)>, \quad \forall u \in {\cal U}, \quad a \in {\cal A}.
\end{eqnarray}

The pairing can be extended to tensor product algebras by setting
\begin{eqnarray} \label{39}
<u \otimes v, a \otimes b> &=& q^{p(v)p(a)} <u,a> <v,b>.
\end{eqnarray}

It is enough to define the pairing between the generating elements of the two algebras. Pairing for any other elements of ${\cal U}$ and ${\cal A}$ follows from (\ref{37}) and the bilinear form inherited by the tensor product.

As a Hopf algebra ${\cal A}$ is generated by the elements $x$, $y$, $\theta$ and a basis is given by all monomials of the form
$$f = x^k y^l \theta^m$$
where $k\in{\mathbb N}_0$, $l,m\in \{0,1,2\}$. Let us denote the dual algebra by ${\cal U}_q$ and its generating elements by $X$, $Y$ and $\Theta$ which are of degree 0,2,1, respectively.

\begin{thm} 
The commutation relations between the generators of the algebra ${\cal U}_q$ dual to ${\cal A}$ are as follows
\begin{eqnarray} \label{40}
[X,Y] = Y, \quad [Y, \Theta] = 0, \quad [\Theta, X] = 2 \,\Theta, \qquad Y^3=0=\Theta^3.
\end{eqnarray}
\end{thm}

\begin{pf}
The pairing is defined through the tangent vectors as follows
\begin{eqnarray} \label{41}
<X, f> &=& k \delta_{n,0}\delta_{m,0}, \nonumber\\
<Y, f> &=& \delta_{l,1}\delta_{m,0}, \nonumber\\
<\Theta, f> &=& \delta_{l,0}\delta_{m,1}.
\end{eqnarray}
We also have
\begin{eqnarray} \label{42}
<1_{\cal U}, f> &=& \epsilon_{\cal A}(f) = \delta_{n,0}.
\end{eqnarray}
Using the defining relations one gets
\begin{eqnarray} \label{43}
<XY, f> = (k + 1) \delta_{l,1}\delta_{m,0} \quad \mbox{and} \quad <YX, f> = k \delta_{l,1}\delta_{m,0}
\end{eqnarray}
where differentiation is from the right as this is most suitable for differentiation in this basis. Thus one obtains
\begin{eqnarray*}
<XY-YX, f> = \delta_{l,1}\delta_{m,0} = <Y, f>.
\end{eqnarray*}
The other relations can be obtained similarly.
\end{pf}

\begin{thm} 
The Hopf algebra structure of the algebra ${\cal U}_q$ is given by:\\
(1) the comultiplication $\Delta_{\cal U}:{\cal U}_q\longrightarrow {\cal U}_q\otimes{\cal U}_q$,
\begin{eqnarray} \label{44}
\Delta_{\cal U}(X) &=& X \otimes 1_{\cal U} + 1_{\cal U} \otimes X, \nonumber\\
\Delta_{\cal U}(Y) &=& Y \otimes q^{2X} + 1_{\cal U} \otimes Y, \nonumber\\
\Delta_{\cal U}(\Theta) &=& \Theta \otimes q^{2X} + 1_{\cal U} \otimes \Theta.
\end{eqnarray}
(2) the counity $\epsilon_{\cal U}:{\cal U}_q\longrightarrow {\mathbb C}$,
\begin{eqnarray} \label{45}
\epsilon_{\cal U}(X) = 0, \quad \epsilon_{\cal U}(Y) = 0, \quad \epsilon_{\cal U}(\Theta) = 0.
\end{eqnarray}
(3) the coinverse $\kappa_{\cal U}:{\cal U}_q\longrightarrow {\cal U}_q$,
\begin{eqnarray} \label{46}
\kappa_{\cal U}(X) = - X \quad \kappa_{\cal U}(Y) = - Y q^X, \quad \kappa_{\cal U}(\Theta) = - \Theta q^X.
\end{eqnarray}
\end{thm}

\begin{pf}
Using the rule in (\ref{3}), it is easily shown that the operator $\Delta_{\cal U}$ (and others) leave invariant the relations (\ref{40}).
We will only obtain the actions of $\Delta_{\cal U}$, $\epsilon_{\cal U}$ and $\kappa_{\cal U}$ on $Y$. The others can be obtained similarly. So, we assume that the action of $\Delta_{\cal U}$ on $Y$ is
$\Delta_{\cal U}(Y)= Y \otimes K_1 + K_2 \otimes Y$. Then the commutation relations in ${\cal A}$ (or comultiplication in ${\cal U}$) will imply that $K_1=q^2K_2$. Indeed,
\begin{eqnarray*}
<\Delta_{\cal U}(Y), x\otimes y> &=& <Y \otimes K_1 + K_2 \otimes Y, x\otimes y>\\
&=& <Y,x><K_1,y> + <K_2,x><Y,y> \\
&=& <K_2,x>
\end{eqnarray*}
and
\begin{eqnarray*}
<\Delta_{\cal U}(Y), y\otimes x> &=& <Y \otimes K_1 + K_2 \otimes Y, y\otimes x>\\
&=& <Y,y><K_1,x> + q^2 <K_2,y><Y,x> \\
&=& <K_1,x>.
\end{eqnarray*}
But $xy=qyx$, therefore $K_1=q^2K_2$. Finally, considering $<Y,x^ky>$ and $<Y,yx^k>$ and taking $K_2=1_{\cal U}$ one can find the action of $\Delta_{\cal U}$ on $Y$.

The action of $\epsilon_{\cal U}$ on $Y$ is
\begin{eqnarray*} 
\mu \circ (\epsilon_{\cal U} \otimes 1_{\cal U})\Delta_{\cal U}(Y) = 1_{\cal U}(Y) = Y
&\Longrightarrow & \epsilon_{\cal U}(Y)q^{2X} + Y = Y\\
&\Longrightarrow & \epsilon_{\cal U}(Y) = 0.
\end{eqnarray*}

The action of $\kappa_{\cal U}$ on $Y$ is
\begin{eqnarray*} 
\mu \circ (\kappa_{\cal U} \otimes 1_{\cal U})\Delta_{\cal U}(Y) = \epsilon_{\cal U}(Y) = 0
&\Longrightarrow & \kappa_{\cal U}(Y)q^{2X} + Y = 0\\
&\Longrightarrow & \kappa_{\cal U}(Y) = - Y q^X.
\end{eqnarray*}
\end{pf}

We can now transform this algebra to the form obtained in Section 4 by having the following definitions:
\begin{eqnarray*}
T_x = \frac{q^X-1_{\cal U}}{1-q^2}, \quad T_y = Y, \quad T_\theta = q^{2X}\Theta
\end{eqnarray*}
which are consistent with the commutation relations and the Hopf structures.

\section{The Z$_3$-graded quantum supergroup GL$_q(2/1)$}\label{sec6}

In this section, we will consider a quantum matrix $T$ in Z$_3$-graded quantum superspace and we will obtain commutation relations between the matrix elements of $T$.

Let $T$ be a 3x3 matrix in Z$_3$-graded superspace,
\begin{eqnarray} \label{47}
T = \begin{pmatrix} a & \gamma_1 & \beta_1 \\ \beta_2 & b & \gamma_2 \\ \gamma_3 & \beta_3 & c \end{pmatrix} =(t_{ij})
\end{eqnarray}
where the generators $a$, $b$, $c$ are of degree 0, the generators $\beta_i$ and $\gamma_i$ are of degree 1 and 2, respectively. We will denote by $GL(2|1)$ the group of such matrices.

We now consider the left coaction
\begin{equation} \label{48}
\delta_L: {\mathbb R}_q(2|1) \longrightarrow GL(2|1) \otimes {\mathbb R}_q(2|1), \quad \delta_L(X_i) = t_{ik} \otimes X_k
\end{equation}
(sum over repeated indices). This coaction is an algebra homomorphism and the elements $\delta_L(X_i)$ should satisfy the relations (\ref{1}).
The action of $\delta_L$ on the dual generators is as follows
\begin{equation} \label{49}
\delta_L(\hat{X}_i) = (\tau\otimes {\sf d})(t_{ik} \otimes X_k) = q^{p(t_{ik})} \, t_{ik} \otimes \hat{X}_k
\end{equation}
(sum over repeated indices) where $\hat{X}=(\hat{X}_i)=(\xi, \eta, z)^t$. These elements $\delta_L(\hat{X}_i)$ must satisfy the relations (\ref{2}).

Consequently, we have the following commutation relations between the matrix elements of $T$:
\begin{eqnarray} \label{50}
a \gamma_1 &=& q\gamma_1 a, \quad a\gamma_2=q\gamma_2 a, \quad a \gamma_3 = q \gamma_3 a, \quad a \beta_1 = \beta_1 a, \quad a \beta_2 = q\beta_2 a, \nonumber\\
a \beta_3 &=& q^2 \beta_3 a, \quad ab = ba + (q^2-1) \gamma_1\beta_2, \quad ac = ca + (q-1) \beta_1\gamma_3,\nonumber\\
b \gamma_1 &=& \gamma_1 b, \quad b\gamma_2=\gamma_2 b, \quad b \gamma_3 = q \gamma_3 b + (1-q^2) \beta_2\beta_3, \\
b \beta_1 &=& q^2 \beta_1 b + (1-q^2) \gamma_1\gamma_2, \quad b \beta_2 = q\beta_2 b, \quad b \beta_3 = q\beta_3 b, \nonumber\\
bc &=& cb + (q-1) \gamma_2\beta_3, \quad c \gamma_1 = q^2 \gamma_1 c + (q^2-1) \beta_1\beta_2, \quad c\gamma_2=\gamma_2 c, \nonumber\\
c \gamma_3 &=& q^2 \gamma_3 c, \quad c \beta_1 = \beta_1 c, \quad c \beta_2 = \beta_2 c + (1-q) \gamma_2\gamma_3, \quad c \beta_3 = q^2 \beta_3,  \nonumber\\
\beta_1\beta_2 &=& \beta_2\beta_1 + (q-q^2) a\gamma_2, \quad \beta_1\beta_3 = q \beta_3\beta_1, \quad \beta_2\beta_3 = q \beta_3\beta_2, \nonumber\\
\gamma_1\gamma_2 &=& q^2 \gamma_2\gamma_1, \quad \gamma_1\gamma_3 = \gamma_3\gamma_1 + (1-q^2) a\beta_3, \quad \gamma_2\gamma_3 = q\gamma_3\gamma_2, \nonumber\\
\gamma_1\beta_1 &=& q^2 \beta_1\gamma_1, \quad \gamma_1\beta_2 = q^2 \beta_2\gamma_1, \quad \gamma_1\beta_3 = q^2 \beta_3\gamma_1, \quad \gamma_2\beta_1 = q \beta_1\gamma_2,  \nonumber\\
\gamma_2\beta_2 &=& q \beta_2\gamma_2, \,\, \gamma_2\beta_3 = \beta_3\gamma_2, \,\, \gamma_3\beta_1 = \beta_1\gamma_3, \quad
\gamma_3\beta_2 = q \beta_2\gamma_3, \quad \gamma_3\beta_3 = \beta_3\gamma_3 \nonumber\\
\beta_1^3 &=& \beta_2^3 = \gamma_2^3 = \gamma_3^3 = 0. \nonumber
\end{eqnarray}
These relations describe a deformation of the algebra of functions on GL$(2|1)$. This algebra is denoted by $Fun(\mbox{GL}_q(2|1))$.

The algebra $Fun(\mbox{GL}_q(2|1))$ is a Z$_3$-graded Hopf algebra. Together with this fact, the properties of the Z$_3$-graded quantum supergroup GL$_q(2|1)$  will be considered in another paper.

\section{An $R$-matrix Formalism}\label{sec7}

Using the relations (\ref{14}), we can find an $R$-matrix satisfying the Z$_3$-graded Yang-Baxter equation. We assume that an $R$-matrix is associated with the quantum Z$_3$-graded superspace ${\mathbb R}_q(2|1)$. Then, we can express the commutation relations, (\ref{14}), between the coordinate functions and their differentials in the form
\begin{eqnarray} \label{52}
q^{p(X_i)} X_i \, {\sf d}X_j = q\hat{R}^{ij}_{kl} \, {\sf d}X_k \, X_l,
\end{eqnarray}
where $\hat{R}=PR$ with the permutation matrix $P$. If we compare it with the relations (\ref{14}), for example, we see that, except the zero entries, $\hat{R}^{21}_{12}=q$, $\hat{R}^{21}_{21}=1-q^2$, for $i=2$, $j=1$. When we find all entries, it takes the following form
\begin{eqnarray} \label{53}
\hat{R} = \begin{pmatrix}
1 & 0 & 0 & 0 & 0 & 0 & 0 & 0 & 0 \\
0 & 0 & 0 & q & 0 & 0 & 0 & 0 & 0 \\
0 & 0 & 1-q^2 & 0 & 0 & 0 & 1 & 0 & 0 \\
0 & q & 0 & 1-q^2 & 0 & 0 & 0 & 0 & 0 \\
0 & 0 & 0 & 0 & 1 & 0 & 0 & 0 & 0 \\
0 & 0 & 0 & 0 & 0 & 1-q^2 & 0 & q & 0 \\
0 & 0 & q^2 & 0 & 0 & 0 & 0 & 0 & 0 \\
0 & 0 & 0 & 0 & 0 & q & 0 & 0 & 0 \\
0 & 0 & 0 & 0 & 0 & 0 & 0 & 0 & 1 \\
\end{pmatrix}=(\hat{R}_{ij}).
\end{eqnarray}

The commutation rules of the generators, $X=(X_i)$, of function algebra on the quantum superspace ${\mathbb R}_q(2|1)$, with an $R$-matrix, can be expressed as follows
\begin{eqnarray} \label{54}
X\otimes X = \hat{R} \, X\otimes X.
\end{eqnarray}

Using the $\hat{R}$-matrix, we can rewrite the relations (\ref{25}) and (\ref{27}) as follows
\begin{eqnarray} \label{55}
\partial_i \, X_j = \delta_{ij} \, + q\hat{R}_{il}^{jk} \, X_l \, \partial_k, \quad
\partial_i \, \partial_j = \hat{R}_{ji}^{lk} \, \partial_k \, \partial_l.
\end{eqnarray}
All relations in Section 3 can be rewritten in terms of the $\hat{R}$-matrix.

\section*{Acknowledgements}
This work was supported in part by {\bf TUBITAK} the Turkish Scientific and Technical Research Council.

\end{document}